\documentclass[onefignum,onetabnum]{siamart190516}

\usepackage{amsmath,amssymb,ifthen,url,graphicx,enumerate,array,theorem,booktabs, tikz}

\newcommand{\Z}{\mathbb{Z}}
\newcommand{\R}{\mathbb{R}}

\newcommand{\IP}{\operatorname{IP}}

\newcommand{\LP}{\operatorname{LP}}

\newcommand{\ceil}[1]{\lceil#1\rceil}
\newcommand{\prt}{\operatorname{Pr_{\emph{t}}}}
\newcommand{\intr}{\operatorname{int}}

\newcommand{\supp}{\operatorname{supp}}
\newcommand{\cone}{\operatorname{cone}}

\newcommand{\rank}{\operatorname{rank}}

\newcommand{\vol}{\operatorname{vol}}
\newcommand{\ba}{\mathbf{a}}
\newcommand{\bb}{\mathbf{b}}
\newcommand{\bc}{\mathbf{c}}
\newcommand{\bd}{\mathbf{d}}
\newcommand{\bg}{\mathbf{g}}
\newcommand{\bh}{\mathbf{h}}

\newcommand{\bu}{\mathbf{u}}
\newcommand{\bv}{\mathbf{v}}
\newcommand{\bw}{\mathbf{w}}
\newcommand{\bx}{\mathbf{x}}
\newcommand{\by}{\mathbf{y}}
\newcommand{\bz}{\mathbf{z}}
\newsiamremark{remark}{Remark}

\newsiamthm{claim}{Claim}

\headers{The distributions of functions related to IP}{T. Oertel, J. Paat, and R. Weismantel}

\title{The distributions of functions related to parametric  integer optimization
}

\author{Timm Oertel\thanks{School of Mathematics, Cardiff University, United Kingdom 
  (\email{oertelt@cardiff.ac.uk}).}
\and Joseph Paat\thanks{Sauder School of Business, University of British Columbia, Canada
  (\email{joseph.paat@sauder.ubc.ca}).}
\and Robert Weismantel\thanks{Department of Mathematics, Institute for Operations Research, ETH Z\"urich, Switzerland 
  (\email{robert.weismantel@ifor.math.ethz.ch}).}}

\usepackage{amsopn}

\begin{document}

\maketitle

\begin{abstract}
We create a framework for studying the asymptotic distributions of functions related to integer linear optimization.
Each of these functions is defined for a fixed constraint matrix and objective vector while the right hand side is treated as input. 
We provide a spectrum of probability-like results that govern the overall asymptotic distribution of a function.
We then apply this framework to the $\IP$ sparsity function, which measures the minimal support of optimal $\IP$ solutions, and the $\IP$ to $\LP$ distance function, which measures the distance between optimal $\IP$ and $\LP$ solutions.
There has been a significant amount of research regarding the extreme values that these functions can attain.
However, less is known about their typical values.
Our results show that the typical values are smaller than the known worst case bounds.
\end{abstract}



\section{Introduction}

%
Let $A \in \Z^{m \times n}$ with $\rank(A) = m$ and $\bc \in \mathbb{Q}^n$ satisfy $\bc^\intercal \bx \le 0$ for all $\bx \in \R^n_{\ge 0} $ such that $A\bx = \mathbf{0} $.
We consider $A$ and $\bc$ to be fixed throughout the paper.
For every $\bb \in \Z^m$, define the integer program
\begin{equation*}
\max\{\bc^\intercal \bz : A\bz = \bb \text{ and } \bz \in \Z^n_{\ge 0}\}.\tag*{IP({\bf b})}
\end{equation*}
The study of $\IP(\bb)$ as $\bb$ varies is referred to as parametric integer programming. %
See Papadimitriou~\cite{P1981} or Eisenbrand and Shmonin~\cite{ES2008}.
The motivation of this paper is to understand $\IP(\bb)$ by studying functions $f$ whose input is $\IP(\bb)$, or equivalently, whose input is a vector $\bb \in \Z^m$.
Such functions include the \emph{integrality gap function}~\cite{AHO2019,DF1989,HS2007}, the \emph{optimal value function}~\cite{G1965,W1981}, the \emph{running time of an algorithm} as a function of $\bb$ \cite{MR2974303,paat2019integrality}, and the \emph{flatness value} \cite{BLPS1999,GC2016}. 
Other examples include the \emph{sparsity function} and the \emph{$\IP$ to $\LP$ distance function}.
Each of the previous functions, when properly normalized, fit into the framework described in this paper.
These functions are well studied in terms of the worst case, e.g., their maximum values.
However, little is known about their distributions, e.g., expected values or how often the worst case occurs.
%
We believe that studying these distributions may lead to improvements in dynamic programs for parametric integer programming, say in the average case.

Let $f:\Z^m \to \R_{\ge 0}\cup\{\infty\}$.
We make the natural assumption that
\begin{equation}\label{eqFeasFin}
f(\bb) < \infty \text{ if and only if } \IP(\bb) \text{ is feasible}.
\end{equation}
In light of the assumption on $A$ and $\bc$ made in the beginning, we see that if $\IP(\bb)$ is feasible, then there exists an optimal solution.
Some choices of $f$ are known to have asymptotically periodic distributions.
Examples include the optimal value function~\cite{G1965} and the sparsity function~\cite{ADOO2017}. 
Underlying the proofs of periodicity is the idea that these functions are well behaved on a family of lattices. 
By exploring these lattice structures in more detail, we can quantify the occurrences of \emph{common values of $f(\bb)$}.
The goal of this paper is to provide lower bounds for these common values.

We quantify common values of $f(\bb)$ using lower asymptotic densities.
For $t \in \Z_{\ge 1}$ and $E \subseteq \Z^m$, define
\[
\prt(E) := \frac{| \{\bb \in E : \|\bb\|_{\infty} \le t \text{ and } f(\bb) < \infty\}  |}
{| \{\bb \in\Z^m: \|\bb\|_{\infty} \le t \text{ and }  f(\bb) < \infty\} |}.
\]
The value $\prt(E)$ is the probability of randomly selecting an integer program $\IP(\bb)$ with $\bb\in E$ among the feasible integer programs with $\bb \in \{-t, \dotsc, t\}^m$.
The \emph{lower asymptotic density of $E$} is
\[
\Pr(E) := \liminf_{t \to \infty} ~ \prt(E).
\]
The value $\Pr(E)$ is the chance of randomly selecting $\IP(\bb)$ with $\bb\in E$ among all feasible integer programs.
The term density is adopted from number theory, see~\cite[Page xii and \S 16]{N2000}.
We use the term density rather than probability because $\Pr(\cdot)$ is not necessarily a probability measure. 
Indeed, it satisfies $\Pr(E) \in [0,1]$ and $\Pr(F) \le \Pr(E)$ if $F \subseteq E$, but not necessarily $\Pr(E \cap F) + \Pr(E \cup F) = \Pr(E) + \Pr(F)$.
We choose to define $\Pr(E)$ as a lower density so that it is well defined for general $f$ and $E$.
However, every limit inferior that we compute is actually a limit. Thus, we often replace  `$\liminf$' by `$\lim$'. 
We are interested in densities of the form
\[
\Pr(f \le \alpha) := \Pr(\{\bb \in \Z^m : f(\bb) \le \alpha\}),
\]
where $\alpha\in\R_{\ge 0}$.
Our first main contribution is Theorem~\ref{thmMain}, which is a set of conditions to bound $\Pr(f \le \alpha)$ for general functions $f$ and values $\alpha$.
The formal result and the intuition behind our proof are presented in Section~\ref{secGeneral} because they require some preliminaries. 
The bounds in Theorem~\ref{thmMain} are in terms of $m$ and the determinants of the submatrices of $A$. 
We denote the largest absolute value of these determinants and their greatest common divisor by
\begin{equation}\label{eqGCD}
\begin{array}{rcl}
\delta & := & \max ~  \{|\det(B)|: B \text{ is an } m \times m \text{ submatrix of } A \}~~ \text{and}\\[.15 cm]
 \gamma  &:=& \gcd~( \{|\det(B)|: B \text{ is an } m \times m \text{ submatrix of } A\}).
\end{array}
\end{equation}
Our second main contribution is an application of Theorem~\ref{thmMain} to bound the asymptotic densities for the sparsity and distance functions.
%

\subsection{The sparsity function $\sigma$}\label{subsecSparsity}

For $\bz \in \R^n_{\ge 0}$, set $\supp(\bz) := \{i \in \{1, \dotsc, n\} : \bz_i > 0\}$.
The \emph{minimum sparsity of an optimal solution to $\IP(\bb)$} is
\[
\sigma(\bb) := \min\{|\supp(\bz)| : \bz \text{ is an optimal feasible solution to } \IP(\bb)\}.
\]
If $\IP(\bb)$ is infeasible, then $\sigma(\bb) := \infty$.
The function $\sigma$ has been used to measure distance between linear codes~\cite{APY2009,A1997} and sparsity in combinatorial problems~\cite{CCD2007,KK1982}.

It was shown by Aliev et al.~\cite{ADEOW2018, ADOO2017} that if $\sigma(\bb) < \infty$, then 
\begin{equation}\label{eqSuppUpperBound}
\sigma(\bb) \le m + \log_{2}(\gamma^{-1}\cdot \sqrt{\det(A{\displaystyle{A^\intercal}})}) \le 2m\log_2(2\sqrt{m}\cdot\|A\|_{\infty}),
\end{equation}
where $\|A\|_{\infty}$ denotes the largest absolute entry of $A$.
See also Eisenbrand and Shmonin~\cite{ES2006}.
In general, there is not much room to improve~\eqref{eqSuppUpperBound}.
For any $\epsilon > 0$, Aliev et al.~\cite{ADEOW2018} provide an example of $A$ and $\bb$ such that
\[
\sigma(\bb) \ge  m \log_2(\|A\|_{\infty})^{1/(1+\epsilon)}.
\]

If $\bc = \mathbf{0}^n$, then $\sigma(\bb)$ quantifies the sparsest \emph{feasible} solution to $\IP(\bb)$. 
Upper bounds on $\sigma(\bb)$ under this assumption were studied in \cite{AlAvDeOe19,ADOO2017}.
Furthermore, Oertel et al.~\cite{OPW2019} showed that asymptotic densities of $\sigma$ can be bounded using the minimum absolute determinant of $A$ or the `number of prime factors' of the determinants.
If, in addition, $A$ has the Hilbert basis property (i.e., if the columns of $A$ correspond to a Hilbert basis of the cone generated by $A$), then bounds on $\sigma(\bb)$ can be given solely in terms of $m$.
Cook et al.~\cite{CookFS1986} showed that if $\sigma(\bb) < \infty$, then $\sigma(\bb) \le 2m-1$; this was improved to $\sigma(\bb) \le 2m-2$ by Seb\H{o}~\cite{Sebo1990}.
Bruns and Gubeladze proved that $\Pr(\sigma \le 2m-3) = 1$~\cite{BG2004}, and Bruns et al.~\cite{BrunsGHMW1999} gave an example such that $\sigma(\bb) \ge (7/6) m $.

We show that $\sigma(\bb)$ is often smaller than the best known worst case bound~\eqref{eqSuppUpperBound}.
\begin{theorem}\label{thmSuppProb}
For each $k \in \{0, \dotsc,  \ceil{\log_2(\gamma^{-1} \cdot \delta)}\}$, it holds that
\[
\Pr\left( \sigma \le m + k \right) 
\ge  \min\bigg\{1, ~\frac{2^k}{\gamma^{-1} \cdot \delta}\bigg\}.
\]
In particular, $ \Pr\left(\sigma \le m + \log_2(\gamma^{-1} \cdot \delta)\right) = 1$.
\end{theorem}

The Cauchy-Binet formula (see~\cite[Section 0.8.7]{HJ2012}) shows that $\delta \le \sqrt{\det(A{\displaystyle{A^\intercal}})}$, and the inequality is strict if $A$ has at least two invertible submatrices.
Hence, the density bounds in Theorem~\ref{thmSuppProb} are often smaller than the worst case bound~\eqref{eqSuppUpperBound}.
Our result can be refined when $\bc = \mathbf{0}^n$.
See Remark~\ref{remSparsFeas}.

\subsection{The distance function $\pi$}\label{subsecProx}

The $\IP$ to $\LP$ distance function measures the distance between optimal solutions to $\IP(\bb)$ and optimal solutions to its linear relaxation
\begin{equation*}
\max\{ \bc^\intercal \bx : A\bx = \bb \text { and } \bx \in \R^n_{\ge 0}\}. \tag*{LP({\bf b})}
\end{equation*}
Whenever we consider $\IP$ to $\LP$ distance we assume, for ease of presentation, that the optimal solution to $\LP(\bb)$ is unique for all feasible $\bb$.
Note that this can always be achieved by perturbing $\bc$;
see Remark~\ref{remark:uniqueOptima} for more on this assumption and its implications.
Let $\bx^*(\bb)$ denote the unique optimal solution to $\LP(\bb)$.
Define the distance function to be
\[
 \pi(\bb)
 :=  \min \left\{\|\bx^*(\bb) - \bz^*\|_1 : \bz^* \text{ is an optimal solution to } \IP(\bb) \right\}.
\]
If $\IP(\bb)$ is infeasible, then $\pi(\bb) := \infty$.

The distance between solutions to $\IP(\bb)$ and $\LP(\bb)$ is a classic question in IP theory that has been used to measure the sensitivity of optimal $\IP$ solutions~\cite{BJ1977, BJ1979,CGST1986} and to create efficient dynamic programming algorithms~\cite{EW2018,JR2018}.
Eisenbrand and Weismantel~\cite{EW2018} showed that if $\pi(\bb) < \infty$, then $\pi(\bb) \le m(2m\|A\|_{\infty}+1)^m$.
By modifying their proof\footnote{The proof of~\eqref{eqProxUpperBound} is the same as~\cite[Theorem 3.1]{EW2018} except the $\|\cdot\|_{\infty}$-norm is replaced by the norm $\|\bx\|_{*} := \|B^{-1} \bx\|_{\infty}$, where $B$ is an $m\times m$ submatrix of $A$ satisfying $|\det(B)| = \delta$.}, it can be shown that if $\pi(\bb) < \infty$, then
\begin{equation}\label{eqProxUpperBound}
\pi(\bb) \le m (2m+2)^m \delta.
\end{equation}
See~\cite{AHO2019,BJ1977,BJ1979,CGST1986,PWW2018,LX2019} for other bounds on $\pi$.
It is not known if the bound in~\eqref{eqProxUpperBound} is tight.
In the case $m=1$, Aliev et al.~\cite{AHO2019} provide a tight upper bound on the related distance function
\[
 \pi^{\infty}(\bb)
 :=  \min \left\{\|\bx^*(\bb) - \bz^*\|_{\infty} : \bz^* \text{ is an optimal solution to } \IP(\bb) \right\}.
\]

Gomory proved that the value function of $\IP(\bb)$ is asymptotically periodic~\cite{G1965}, see also Wolsey~\cite{W1981}.
Using his results along with Theorem~\ref{thmMain}, one can prove that $\Pr(\pi \le (m+1)\gamma^{-1} \cdot \delta) = 1$.
We provide a refined density analysis in Theorem~\ref{thmMainProx} \emph{(a)}. 
Theorem~\ref{thmMainProx} \emph{(b)} bounds densities in terms of $ \pi^{\infty}$.
\begin{theorem}\label{thmMainProx} For each $k  \in \{0, \dotsc, \gamma^{-1} \cdot \delta - 1\}$, it holds that
\begin{enumerate}[(a)]
\item  $\displaystyle\Pr\left(\pi \le m  \gamma^{-1} \cdot \delta \cdot \frac{k}{k+1} +k \right)  \ge  \frac{k+1}{\gamma^{-1} \cdot \delta}$  and\\[.25 cm]
\item $\displaystyle  \Pr\left(\pi^{\infty} \le  \gamma^{-1} \cdot \delta \cdot \frac{k}{k+1} \right) \ge  \frac{k+1}{\gamma^{-1} \cdot \delta}$.\\[.1 cm]
\end{enumerate}
In particular, $\Pr(\pi \le (m+1) (\gamma^{-1} \cdot \delta-1))=1$ and $\Pr\left(\pi^{\infty} \le  \gamma^{-1} \cdot \delta -1\right)  = 1$.
\end{theorem}

Theorem~\ref{thmMainProx} \emph{(b)} partially resolves Conjecture 1 in~\cite{PWW2018}, which states that $\pi^{\infty}$ can be bounded in terms of the largest minor of $A$ and independently of the number of constraints $m$ and the dimension $n$.
Together with Hadamard's inequality (see, e.g.,~\cite[Corollary 7.8.3]{HJ2012}), Theorem~\ref{thmMainProx} can be used to bound the typical distance between solutions to $\IP(\bb)$ and $\LP(\bb)$ in terms of $\|A\|_{\infty}$ rather than $\delta$.
\begin{corollary}
The function $\pi$ satisfies 
\[
\Pr(\pi \le (m+1) \cdot ( \sqrt{m} \|A\|_{\infty})^m) = 1.
\]
\end{corollary}
%

\subsection{Outline and notation}\label{subsecOutline}

%
Section~\ref{secGeneral} provides a general framework for upper bounding $\Pr(f \le \alpha)$ and proves the fundamental Theorem~\ref{thmMain}.
Preliminaries about optimal solutions to $\IP(\bb)$ are given in Section~\ref{secOptimal}.
We use these preliminaries in Section~\ref{secSupp} to prove Theorems~\ref{thmSuppProb} and~\ref{thmMainProx}.

We view $A $ as a matrix and as a set of column vectors in $\Z^m$, so $B \subseteq A$ means $B$ is a subset of the columns of $A$.
For $K \subseteq \R^m$ and $\bd \in \R^m$, define $K + \bd :=\{ \bb + \bd : \bb \in K\}$.
The $k$-dimensional vector of all zeros is denoted by $\mathbf{0}^k$, and the vector of all ones is denoted by $\mathbf{1}^k$.
When multiplying a matrix $B \subseteq \Z^{m}$ and a vector $\by \in \R^{B}$ as $B\by$, we use $\by_{\bb}$ to denote the component of $\by$ corresponding to $\bb\in B$. 
For $P \subseteq \R^m$, we use $\cone(P)$ to denote the \emph{convex cone} generated by $P$ and $\intr(P)$ to denote the \emph{interior of $P$}.
The \emph{dimension of $P$} is the dimension of the affine hull of $P$.

A set $\Lambda \subseteq \Z^m$ is a \emph{lattice} if $\mathbf{0}^m \in \Lambda$, $\bb + \bd \in \Lambda$ if $\bb,\bd\in \Lambda$, and $-\bb \in \Lambda$ if $\bb \in \Lambda$. 
If $\bb \in \Z^m$ and $\Lambda$ is a lattice, then $\Gamma = \bb + \Lambda$ is an \emph{affine lattice}.
The \emph{dimension of $\Gamma$} is the largest number of linearly independent vectors in $\Lambda$.
The \emph{determinant} of an $m$-dimensional affine lattice $\Gamma$ is $\det(\Gamma) := |\det(B)|$, where $B \in \Z^{m\times m}$ is any matrix such that $\Lambda = B \cdot \Z^m$.
An $m$-dimensional lattice $\Lambda$ induces an \emph{equivalence relationship $\equiv_{\Lambda}$} on $\Z^m$, where $\bb \equiv_{\Lambda} \bd$ if and only if $\bb - \bd \in \Lambda$.
The number of equivalence classes induced by $\equiv_{\Lambda}$ is $\det(\Lambda)$~\cite[Page 22]{GruLek87}.
We refer to~\cite{AS1986} and~\cite[Chapter VII]{barv2002} for more on lattices.

A particular lattice that we use throughout is
\begin{equation}\label{eqALattice}
\Lambda := A \cdot \Z^n.
\end{equation}
Note that $\det(\Lambda) = \gamma$, where $\gamma$ is defined in~\eqref{eqGCD}.
For completeness, we give a short proof.
Let $B\in\Z^{m\times m}$ be such that $\Lambda=B\cdot\Z^m$.
Thus, $|\det(B)| = \det(\Lambda)$.
Let $D$ be any subset of $m$ columns of $A$.
There exists a matrix $U\in\Z^{m\times m}$ such that $D = B U$ because $A\subseteq\Lambda$.
Thus, $\det(B)\mid\det(D)$.
It follows that $\det(B)\mid\gamma$ because $D$ was chosen arbitrarily.
Conversely, there exists a matrix $V\in\Z^{n\times m}$ such that $B=AV$ because $\Lambda=A\cdot\Z^n$.
The Cauchy-Binet formula states that
\[
\det(B) = \sum_{\substack{I \subseteq \{1, \dotsc, n\} \\ |I| = m}} \det(A_I) \cdot \det(V_I),
\]
where $A_I$ and $V_I$ denote the matrices formed by the columns of $A$ and the rows of $V$ indexed by $I$, respectively. 
Thus, $\gamma\mid\det(B)$.
%

\section{Asymptotic densities for general functions}\label{secGeneral}

%
Let $f : \Z^m \to \R_{\ge 0} \cup \{\infty\}$ satisfy~\eqref{eqFeasFin}, $\alpha \in \R$, and $\Lambda = A \cdot \Z^n$.
%
The key idea behind how we lower bound $\Pr(f \le \alpha)$ is to exploit potential local periodic behavior of $f$.
We briefly outline this idea below.
We say that a right hand side $\bb \in \Z^m$ is `good' if $f(\bb) \le \alpha$. 
Assumption~\eqref{eqFeasFin} implies that a good right hand side must be in $\Lambda$, so we may restrict ourselves to consider $\bb$ in $\Lambda$ rather than in $\Z^m$.
%

%
First, we cover $\cone(A)$ by \emph{simplicial cones} $\cone(A^1),$ $\dotsc, \cone(A^s)$, where $A^1,$ $\dotsc, A^s \subseteq A$.
The density of good vectors in $\cone(A)$ is larger than the minimum density of good vectors in any $\cone(A^i)$.
Hence, it suffices to lower bound the density of good vectors in each $\cone(A^i)$ individually.
Not every $\bb \in \cone(A^i) \cap \Lambda$ is feasible, but one can show that there exists a vector $\bd^i \in \cone(A^i) \cap \Z^m$ such that $\IP(\bb)$ is feasible for all $\bb\in[\cone(A^i)+\bd^i]\cap\Lambda$.
This phenomenon relates to the \emph{Frobenius number}, see~\cite{AOW2015,S2012}.
Motivated by these `deep' regions, we use \emph{Ehrhart theory} to show that the density of good vectors in $\cone(A^i) + \bd^i$ is equal to the density of good vectors in $\cone(A^i)$.
See Lemma~\ref{lem:EhrhartTheory}.

Next, we consider the sublattice $\Gamma^i = A^i \cdot \Z^m$, which serves as a natural candidate for quantifying periodicity within $\cone(A^i)$.
The lattice $\Lambda$ is covered by the disjoint affine lattices $\{\Gamma^i +\bg : \bg \in \Lambda/\Gamma^i\}$.
Instead of computing the density of good vectors in $\cone(A^i) + \bd^i$, we count the number of disjoint affine lattices with the property that all vectors in $[\cone(A^i) + \bd^i]\cap[ \Gamma^i+\bg ]$ are good.
See \eqref{thmMain:condition}.

We now formalize the steps above. 
We say that matrices $A^1, \dotsc, A^s \subseteq A$ form a \emph{simplicial covering of $\cone(A)$} if each $A^i$ is invertible, i.e., $\cone(A^i)$ is simplicial, and
\[
\cone(A) = \bigcup_{i=1}^s \cone (A^i).
\]
These coverings always exist due to Carath\'{e}odory's theorem.
The cones in a simplicial covering may overlap nontrivially.
In order to prevent double counting, we triangulate the cones using the next lemma. 
We omit the proof as it follows from standard results on triangulations and subdivisions.
See~\cite[Page 332]{barv2002} or~\cite[Chapter 9]{Z95}.
\begin{lemma}\label{lem:unimodularCovering}
Let $A^1,\ldots,A^s\in\Z^{m\times m}$ be square matrices of rank $m$.
There exist $m$-dimensional rational polyhedral cones $C^1,\ldots, C^\ell\subseteq \R^m$ such that
\begin{enumerate}[(a)]
\smallskip
\item $\bigcup_{i=1}^s \cone(A^i) = \bigcup_{j=1}^\ell C^j$,
\smallskip
\item $\intr(C^j) \cap \intr(C^k)=\emptyset$ for distinct $j,k \in \{1, \dotsc, \ell\}$, and
\smallskip
\item $C^j \subseteq\cone(A^i)$ or $\intr(C^j) \cap \cone (A^i)=\emptyset$ for all $i \in \{1 \dotsc, s\} $ and $j \in \{1, \dotsc, \ell\}$.
\end{enumerate}
\end{lemma}
For functions $g,h:\R_{>0}\to\R_{>0}$, we write 
\[
g\sim h ~~ \text{if}~\lim_{t\to\infty}\frac{g(t)}{h(t)}=1 \qquad\text{and}\qquad g \precsim h ~~\text{if}~\limsup_{t\to\infty}\frac{g(t)}{h(t)}\le1.
\]
For a $q$-dimensional set $P \subseteq \R^m$, we denote the $q$-dimensional Lebesgue measure by $\vol_q(P)$.
The next lemma will enable us to compare densities, and it is a variation of classic results in Ehrhart theory.
See~\cite[Theorem~7]{McMullen78} and~\cite[Theorem~1.2]{HenLin15}.
%
%

%
\begin{lemma}\label{lem:EhrhartTheory}
Let $P \subseteq \R^m$ be a $q$-dimensional rational polytope and $\Gamma\subseteq\Z^m$ an $m$-dimensional affine lattice.
There exists a constant $\eta_{P,\Gamma}>0$ such that 
\[
| t P\cap\Gamma | \precsim \eta_{P,\Gamma} \cdot t^q.
\]
If $q=m$, then $\eta_{P,\Gamma}=\vol_m (P)/{\det(\Gamma)}$ and 
\[
| t P\cap\Gamma| \sim \eta_{P,\Gamma} \cdot t^m.
\]
\end{lemma}

Define the lattices 
\begin{equation}\label{eqGammaLattice}
\Gamma^{i} := A^i \cdot \Z^m \quad \forall ~ i \in \{1, \dotsc, s\}
\end{equation}
with corresponding equivalence relations $\equiv_{\Gamma^i}$.
Observe that $\det(\Gamma^i) = |\det(A^i)|$ and that $\Gamma^i$ is a sublattice of $\Lambda$ for each $i \in \{1, \dotsc, s\}$.
Hence, the relation $\equiv_{\Gamma^i}$ induces a quotient group $\Lambda / \Gamma^i$ with cardinality
\begin{equation}\label{eqNormalizedGCD}
|\Lambda / \Gamma^i| = \det (\Gamma^i)/\det (\Lambda) = \gamma^{-1} \cdot |\det(A^i)|.
\end{equation}
In other words, $\equiv_{\Gamma^i}$ partitions $\Lambda$ into $\gamma^{-1} \cdot |\det(A^i)|$ many different equivalence classes.
We are now prepared to formally state our first main result.

\begin{theorem}\label{thmMain}
Let $f$ satisfy~\eqref{eqFeasFin}, $\alpha \in \R$, and  $A^1, \dotsc, A^s$ be a simplicial covering of $\cone(A)$. 
Set $\Lambda=A\cdot\Z^n$ .
For each $i \in \{1, \dotsc, s\}$, let $\bd^i \in \cone(A^i)\cap\Z^m$, and define $\Gamma^i=A^i\cdot\Z^m$ and
\begin{equation}\label{thmMain:condition}
\beta_i := \left|\left\{\bg \in \Lambda/\Gamma^i : 
\max\left\{f(\bb) :
\hspace{-.15 cm}
\begin{array}{l}\bb \equiv_{\Gamma^i} \bg,\\[.05 cm]
 \bb\in \cone(A^i) + \bd^i
 \end{array}
 \hspace{-.1 cm}
\right\} \le \alpha\right\}\right|.
\end{equation}
It holds that
\begin{equation}\label{thmMainResultOne}
\Pr\left(f \le \alpha\right) ~\ge~ \min_{i=1,\ldots,s} ~ \frac{\beta_i}{\gamma^{-1}\cdot\det(\Gamma^i)}.
\end{equation}
\end{theorem}

\proof
It follows from~\eqref{eqFeasFin} that if $\bb \in \Z^m$ and $f(\bb) < \infty$, then $\bb \in \Lambda \cap \cone(A)$. 
Therefore, 
\[
 \prt(  f \le \alpha)
= 
\frac{|\{\bb\in\Lambda \cap \cone(A) \;:\; \|\bb\|_\infty \le t \text{ and } f(\bb) \le \alpha\}|}{|\{\bb\in\Lambda \cap \cone(A)  \;:\; \|\bb\|_\infty \le t \text{ and } f(\bb)<\infty\}|}
\quad \forall ~ t \in \Z_{\ge 0}.
\]
By Lemma~\ref{lem:unimodularCovering}, we can cover $\cone(A)$ by rational polyhedral cones $C^1, \dotsc, C^{\ell}$ such that $\intr(C^j)\cap\intr(C^k)=\emptyset$ for distinct $j,k\in\{1,\dots,\ell\}$ and either $C^j\subseteq\cone(A^i)$ or $\intr(C^j)\cap\cone(A^i)=\emptyset$ for all $i\in\{1,\ldots,s\}$ and $j\in\{1,\dots,\ell\}$.
For each  $j \in \{1, \dotsc, \ell\}$, define the truncated cone $P^j:= C^j\cap[-1,1]^m$.
By Lemma~\ref{lem:EhrhartTheory}, there exist positive constants $\eta_j$ and $\eta_{jk}$ such that $| \Lambda \cap tP^j  | \sim \eta_j\,t^m $ and $| \Lambda \cap t(P^j \cap P^k)| \precsim \eta_{jk}\,t^{m-1}$ for any intersection $P^j \cap P^k$ satisfying $j\neq k$.
Asymptotic densities are defined through limits. 
Thus, we may neglect any low-dimensional intersections in the covering of $\cone(A)$ by $C^1, \dotsc, C^{\ell}$ and instead treat the covering as a partition.
We have
\begin{equation}\label{eq:mainThm:Reduction}
\begin{aligned}
& \Pr(f \le \alpha) = \lim_{t\to\infty} \prt(f \le \alpha)\\
= &  \lim_{t\to\infty} \sum_{j=1}^\ell ~ \frac{~~~~~~~|\{\bb\in \Lambda \cap tP^j: f(\bb) \le \alpha\}|}{\sum_{k=1}^{\ell}|\{\bb\in\Lambda \cap tP^k :  f(\bb)<\infty\}|} \\[.125 cm]
\ge &  \lim_{t\to\infty} \sum_{j=1}^\ell ~ \frac{~~~~~~~\,|\{\bb\in \Lambda \cap tP^j: f(\bb) < \infty\}|}{\sum_{k=1}^{\ell}|\{\bb\in\Lambda \cap tP^k  :  f(\bb)<\infty\}|} \cdot \frac{|\{\bb\in \Lambda \cap tP^j: f(\bb) \le \alpha\}|}{|\Lambda \cap tP^j |} \\[.125 cm]
\ge &  \lim_{t\to\infty} ~~ \min_{j = 1, \dotsc, \ell} \frac{|\{\bb\in \Lambda \cap tP^j:  f(\bb) \le \alpha\}|}{|\Lambda \cap tP^j  |} ,\\[.125 cm]
= &  \min_{j = 1, \dotsc, \ell} ~~\lim_{t\to\infty}  \frac{|\{\bb\in \Lambda \cap tP^j:  f(\bb) \le \alpha\}|}{|\Lambda \cap tP^j  |}.
\end{aligned}
\end{equation}
The second equation in~\eqref{eq:mainThm:Reduction} follows because $C^1, \dotsc, C^{\ell}$ partition $\cone(A)$.
The first inequality in~\eqref{eq:mainThm:Reduction} follows because $\{\bb\in\Lambda \cap tP^j  : ~ f(\bb)<\infty\}$ is a subset of $\Lambda \cap tP^j$; thus, it has a smaller cardinality.
The final equation in~\eqref{eq:mainThm:Reduction} holds because the minimum is taken over a finite index set.

Consider a cone $C^j$, where $j\in\{1,\ldots, \ell\}$.
There exists an $i\in\{1,\ldots,s\}$ such that $C^j \subseteq\cone(A^i)$.
In what remains, we prove that
\begin{equation}\label{eqlastStep}
\lim_{t\to\infty} \frac{|\{\bb\in \Lambda \cap  tP^j \;:\; f(\bb) \le \alpha\}|}{| \Lambda \cap tP^j |} \ge \frac{\beta_i}{\gamma^{-1}\cdot \det(\Gamma^i)}.
\end{equation}
The main statement \eqref{thmMainResultOne} follows immediately after combining~\eqref{eq:mainThm:Reduction} and~\eqref{eqlastStep}.

By Lemma~\ref{lem:EhrhartTheory}, the proportion of vectors in $tP^j$ that are also in $\Lambda$ is
\begin{equation}\label{eqAsym0}
| \Lambda \cap tP^j  | \sim t^m \frac{\vol_m(P^j)}{\det (\Lambda)}.
\end{equation}
Similarly, for each $\bg\in\Lambda / \Gamma^i$, the proportion of vectors in $tP^j$ that are in the affine lattice $ \Gamma^i + \bg $ is 
\begin{equation}\label{eqAsym1}
|[ \Gamma^i+\bg ] \cap  tP^j | \sim  t^m \frac{\vol_m(P^j)}{\det (\Gamma^i)}.
\end{equation}
The vectors in $\Gamma^i +\bg$ that are contained in $tP^j \setminus [tP^j+\bd^i]$ lie on a finite number of hyperplanes parallel to the faces of $C^j$.
The number of these hyperplanes is independent of $t$.
%
Thus, by Lemma~\ref{lem:EhrhartTheory},
there exists a constant $\mu>0$ such that
\begin{equation}\label{mainThmTranslatedErhard}
|[\Gamma^i +\bg] \cap [tP^j \setminus [ tP^j + \bd^i] ]| \precsim \mu \cdot t^{m-1}.
\end{equation}
Looking at the difference of \eqref{eqAsym1} and \eqref{mainThmTranslatedErhard}, we obtain 
\begin{equation}\label{eqAsym2}
|  [ \Gamma^i + \bg ] \cap tP^j \cap [tP^j + \bd^i] | \sim  t^m \frac{\vol_m(P^j)}{\det (\Gamma^i)}.
\end{equation}

Set
\[
X^i :=\left\{\bg \in \Lambda/\Gamma^i : 
\max\left\{f(\bb) :
\hspace{-.15 cm}
\begin{array}{l}\bb \equiv_{\Gamma^i} \bg,\\[.05 cm]
 \bb\in \cone(A^i) + \bd^i
 \end{array}
 \hspace{-.1 cm}
\right\} \le \alpha\right\}.
\]
The equation $\beta_i=|X^i|$ holds because of~\eqref{thmMain:condition}.
For each $\bg \in X^i$, it follows that
\[
[ \Gamma^i+\bg ] \cap tP^j
 \supseteq \{\bb\in [\Gamma^i + \bg ] \cap tP^j  : f(\bb) \le \alpha\}
  \supseteq  [ \Gamma^i + \bg ] \cap tP^j \cap [tP^j + \bd^i] .
   \]
   Relations~\eqref{eqAsym1} and~\eqref{eqAsym2} show that the cardinalities of the first and last sets are asymptotically equal. 
Thus,
\begin{equation}\label{eqAsym4}
|\{\bb\in  [\Gamma^i + \bg ] \cap  tP^j : f(\bb) \le \alpha\}| \sim t^m \frac{\vol_m(P^j)}{\det (\Gamma^i)}.
\end{equation}

Every $\bb \in \Lambda\cap tP^j$ belongs to exactly one of the $\gamma^{-1} \cdot \det(\Gamma^i)$ many equivalence classes defined by the relation $\equiv_{\Gamma^i}$.
Therefore,
\begin{equation*} 
|\{\bb\in \Lambda \cap  tP^j  : f(\bb) \le \alpha\}| = \sum_{\bg \in \Lambda / \Gamma^i} |\{\bb\in [\Gamma^i + \bg ] \cap tP^j :f(\bb) \le \alpha\}|.
\end{equation*}
Combining this equation with~\eqref{eqAsym0} and~\eqref{eqAsym4}, we see that
\[
\begin{array}{r@{\hskip .05 cm}rl}
%
\displaystyle \lim_{t\to\infty} \frac{|\{\bb\in \Lambda \cap tP^j : f(\bb) \le \alpha\}|}{| \Lambda \cap tP^j  |}
& = &\displaystyle \lim_{t\to\infty}  \sum_{\bg\in\Lambda / \Gamma^i}\frac{|\{\bb\in [\Gamma^i + \bg ] \cap  tP^j  : f(\bb) \le \alpha\}|}{| \Lambda \cap tP^j |} \\[.625 cm]
& \ge & \displaystyle \lim_{t\to\infty}  \sum_{\bg\in X^i }\frac{|\{\bb\in [ \Gamma^i + \bg] \cap tP^j  : f(\bb) \le \alpha\}|}{| \Lambda \cap tP^j  |} \\[.625 cm]
%
%
& = &\displaystyle  \frac{|X^i|}{|\Lambda / \Gamma^i|} ~=~ \frac{\beta_i}{\gamma^{-1}\cdot\det (\Gamma^i)}.
\end{array}
\]
This proves~\eqref{eqlastStep}.
\endproof

\section{Preliminaries for results on optimal $\IP$ solutions}\label{secOptimal}
%
The density bounds derived in Theorem~\ref{thmMain} depend on the choice of simplicial covering.
We choose a specific covering related to optimal $\LP$ bases in order to prove Theorems~\ref{thmSuppProb} and~\ref{thmMainProx}.
We say that an invertible matrix $B \subseteq A$ is an \emph{optimal $\LP$ basis matrix} if for all $\bb \in \cone(B) \cap \Z^m$ the problem $\LP(\bb)$ has an optimal solution $\bx^*$ satisfying $\{\ba \in A: \bx^*_{\ba} > 0 \} \subseteq B$.
This section collects properties of optimal $\LP$ basis matrices that we will use when applying Theorem~\ref{thmMain} to $\sigma$ and $\pi$.
We begin with a folklore result.
\begin{lemma}\label{lemLPBasisConstant}
The set of all optimal $\LP$ basis matrices defines a simplicial covering of $\cone(A)$.
\end{lemma}

Let $B$ be an optimal basis matrix.
Gomory showed in~\cite[Theorem 2]{G1965} that $\IP(\bb)$ is feasible if $\bb$ is deep inside $\cone(B)$, that is if $\bb$ is in the set\footnote{Gomory defines the set of deep vectors in terms of the distance from $\bb$ to the boundary of $\cone(B)$, and his set contains $D(B)$. 
Our definition of $D(B)$ is chosen to simplify our proofs.}
\begin{equation}\label{eqBSet}
 D(B) := \{
\bb \in \Lambda : 
B^{-1}\bb \ge 3  \delta \cdot  \mathbf{1}^m
\}.
\end{equation}
Furthermore, he showed that there exists an optimal solution $\bz^*$ to $\IP(\bb)$ whose support is contained in $B$ together with few additional non-basic columns $N=A \setminus B$. This fact is shown in Lemma~\ref{lemNonbasic2}.
More precisely, $\bz^* = \bz^B + \bz^N$, where $\{ \ba \in A:\bz^{B}_{\ba} > 0 \} \subseteq B$ and $|\{ \ba \in A:\bz^{N}_{\ba} > 0 \}| < |\det(B)|$. 
Set $ \Gamma := B \cdot \Z^m$. 
Observe that
\[
\bb = A\bz^* = A\bz^B + A\bz^N ~~ \text{and}~~ \{ \ba \in A:\bz^{B}_{\ba} > 0 \} \subseteq B
\]
imply $A \bz^B \equiv_\Gamma \mathbf{0}^m$ and $A \bz^N \equiv_\Gamma \bb$.
Hence, $\bz^{N}$ is the subvector of $\bz^*$ that ensures $A\bz^* \equiv_\Gamma \bb$.
Gomory also argued that $\bz^{N}$ can be chosen to be a minimal subvector with this property.
By minimal, we mean that there does not exist a vector $\overline{\bz}^N \in \Z^n$ satisfying $\mathbf{0}^n \le \overline{\bz}^{N} \lneq \bz^{N}$ and $A \overline{\bz}^{N} \equiv_\Gamma \bb$.
We denote the set of these minimal vectors $\bz^N$ by
\begin{equation}\label{eqDSet}
 N(B) := \left\{
\bz \in \Z^n_{\ge 0} : 
\hspace{-.15 cm}
\begin{array}{l}
\text{there exists}~ \bb \in D(B) \text{ and } \bz^B \in \Z^n_{\ge 0} \text{ such that }\\[.1 cm]
\begin{array}{l@{\hskip .25 cm}l}
(i) &\{\ba \in A:\bz^{B}_{\ba} > 0 \} \subseteq B,\\[.1 cm]
(ii) &\bz^B+ \bz \text{ is an optimal solution to } \IP(\bb),\\[.1 cm]
%
%
(iii) & A \bw \not\equiv_{\Gamma} A\bz ~\text{for all} ~ \mathbf{0}^n \le {\bw} \lneq \bz \\ 
\end{array}
\end{array}
\hspace{-.15 cm} 
\right\}.
\end{equation}

Next, we show that each $\bz \in N(B)$ is not too large and that the coordinates of $A \bz $ in the coordinate space defined by $B$ are not too large either. 
These results only rely on condition \emph{(iii)} in~\eqref{eqDSet}. 

\begin{lemma}\label{lemNonbasic1}
Let $B \subseteq A$ be an optimal $\LP$ basis matrix and $\bz \in \Z^n_{\ge 0}$.
If $A\bw \not\equiv_\Gamma A\bz$ for all $\mathbf{0}^n \le \bw \lneq \bz$, then
\begin{equation}\label{eqBound1}
\|\bz\|_1 < \gamma^{-1}\cdot|\det(B)|
\end{equation}
and
\begin{equation}\label{eqBound2}
\|B^{-1} A \bz\|_{\infty}  \le \|B^{-1}A\|_{\infty} \cdot \|\bz\|_1 < \gamma^{-1}\cdot \delta.
\end{equation}
Consequently, if $\bw \in \Z^n$ and $\bb\in D(B)$ satisfy $\mathbf{0}^n \le \bw \le \bz$ and $A\bw\equiv_{\Gamma}\bb$, then
\begin{equation}\label{eqBound3}
B^{-1}(\bb - A \bw) \ge (3  - \gamma^{-1} )\delta \cdot \mathbf{1}^m \ge \mathbf{0}^m.
\end{equation}
\end{lemma}
\proof
For two vectors $\by, \by'$ satisfying $\mathbf{0}^n \le \by \lneq \by' \le \bz$ we claim that $A \by \not\equiv_{\Gamma} A \by'$. 
Otherwise, we obtain the contradiction $A\bw \equiv_\Gamma A\bz$ and $\mathbf{0}^n \le \bw \lneq \bz$ for the vector $\bw := \bz - \by + \by'$.
Consider any sequence of $\|\bz\|_1+1$ many vectors satisfying $\mathbf{0}^n = \by^1 \lneq \dotsc \lneq \by^{\|\bz\|_1+1} = \bz$.
Each $A \by^i$ is distinct modulo $\Gamma$. 
By~\eqref{eqNormalizedGCD}, there are $\gamma^{-1}\cdot|\det(B)|$ many equivalence classes modulo $\Gamma$.
Hence, $\|\bz\|_1+1 \le  \gamma^{-1}\cdot|\det(B)|$.

Inequality~\eqref{eqBound2} follows from~\eqref{eqBound1} and
\[
 \|B^{-1}A\|_{\infty} \le \frac{\delta}{|\det(B)|}.
 \]
If the latter inequality is false, then there exists $\ba \in A$ and $\bd \in B$ such that $\by := B^{-1}\ba$ and $\by_{\bd} > \delta / |\det(B)|$. 
However,
\(
|\det(B  \cup \{\ba\}\setminus\{\bd\})| = |\by_{\bd}| \cdot |\det(B)| > \delta ,
\)
which contradicts the definition of $\delta $.
\endproof

It is not hard to see that, for every $\bg \in \Lambda / \Gamma$, there exists at least one vector $\bz^{\bg}\in N(B)$ such that $A \bz^{\bg} \equiv_\Gamma \bg$, which also follows from Gomory's work.
The result~\cite[Theorem 2]{G1965} of Gomory can now be stated in terms of $D(B)$ and $N(B)$: If $\bb \in D(B)$, then there exists a vector $\bz \in N(B)$ such that $\bz^B + \bz$ is an optimal solution to $\IP(\bb)$ for some $\bz^B \in \Z^n_{\ge 0}$ satisfying $\{\ba \in A : \bz^B_{\ba}>0\} \subseteq B$.
The following lemma shows a stronger statement: \emph{any} vector $\bz \in N(B)$ can be extended to an optimal solution to $\IP(\bb)$ in this way for any $\bb \in D(B)$ equivalent to $A\bz$.
Furthermore, if $\bz \in N(B)$ and $ \mathbf{0}^n \le {\bw} \le \bz$, then ${\bw} \in N(B)$.

\begin{lemma}\label{lemNonbasic2}
Let $B \subseteq A$ be an optimal $\LP$ basis matrix, $\bz \in N(B)$, and $\bw \in \Z^n$ satisfy $\mathbf{0}^n \le \bw \le \bz$.
For all ${\bb} \in D(B)$ such that $A \bw \equiv_{\Gamma} {\bb}$, there exists an optimal solution to $\IP(\bb)$ of the form $\bw^{B} + \bw$, where $\bw^{B} \in \Z^n_{\ge 0}$ and $\{\ba \in A: \bw^{B}_{\ba} > 0 \} \subseteq B$.
\end{lemma}

\proof
Define $\bw^{B} \in \R^n$ component-wise to be 
\[
\bw^{B}_{\ba} := \begin{cases} 
[B^{-1}(\bb - A\bw)]_{\ba} & \text{if}~ \ba \in B\\[.1 cm]
0 &\text{if}~ \ba \in A\setminus B.
\end{cases}
\]
Note that $\bw^B \in \Z^n$ because $A\bw \equiv_\Gamma \bb$.
Since $\bz \in N(B)$, we may apply Lemma~\ref{lemNonbasic1} to conclude $\|B^{-1}A\|_{\infty} \cdot \|\bz\|_1 < \gamma^{-1} \cdot \delta$. 
Together with $\|\bw \|_1 \le \|\bz\|_1$ this yields
\[
\|B^{-1}A \bw\|_{\infty} \le \|B^{-1}A\|_\infty \cdot\|\bw\|_{1} \le  \|B^{-1}A\|_\infty \cdot\|\bz\|_{1} \le \gamma^{-1} \cdot \delta.
\]
By~\eqref{eqBound3}, $\bw^B$ is nonnegative. 
Thus, $\bw^B+\bw$ is feasible for $\IP(\bb)$.

It remains to show that $\bw^B+\bw$ is optimal for $\IP(\bb)$.
We use an exchange argument to prove this.
The first step is to compare $\bw$ to a vector derived from an optimal solution to $\IP(\bb)$.
There exists an optimal solution $\by^*$ to $\IP(\bb)$ because the problem is feasible and bounded.
Choose $\by \in \Z^n_{\ge 0}$ minimizing $\|\by\|_1$ such that $A \by \equiv_\Gamma \bb$ and $\by \le \by^*$.
The vector $\by$ must satisfy the assumptions in Lemma~\ref{lemNonbasic1}.
Otherwise, $\|\by\|_1$ was not minimized.
Thus,
\[
\|B^{-1}A{\by}\|_{\infty} < \gamma^{-1} \cdot \delta.
\]
Because $A{\bw} \equiv_\Gamma  A{\by}$, there exists a vector $\bu \in \Z^n$ such that
\(
\{\ba \in A: \bu_{\ba} \neq 0 \} \subseteq B
\)
and
\( 
A({\bw} - {\by} + \bu) = \mathbf{0}^m.
\)
Furthermore, 
\begin{equation}\label{eqBoundu}
\|\bu\|_{\infty} 
= \|B^{-1}A({\bw} - {\by})\|_{\infty} 
\le \|B^{-1}A {\bw}\|_{\infty}+ \|B^{-1}A{\by}\|_{\infty} \le 2 \gamma^{-1} \cdot \delta.
\end{equation}
The second step in the exchange argument is to show that
\begin{equation}\label{eqFINALFINALEq1}
 \bc^\intercal(\by^* - \by  + \bu) \le \bc^\intercal \bw^B
\end{equation}
and
\begin{equation}\label{eqFINALFINALEq}
\bc^\intercal(\bw - \by + \bu) = 0.
\end{equation}
The combination of~\eqref{eqFINALFINALEq1} and~\eqref{eqFINALFINALEq} shows that $\bw^B + \bw$ is optimal for $\IP(\bb)$:
\[
\bc^\intercal \by^* = \bc^\intercal(\by^* - \by  + \bu) + \bc^\intercal (\by - \bu) \le \bc^\intercal \bw^B + \bc^\intercal \bw = \bc^\intercal (\bw^B + \bw).
\]
To prove~\eqref{eqFINALFINALEq1}, define $\by^{B} \in \Z^n$ component-wise to be 
\[
\by^{B}_{\ba} := \begin{cases} 
[B^{-1}(\bb - A\by)]_{\ba} & \text{if}~ \ba \in B\\[.1 cm]
0 &\text{if}~ \ba \in A\setminus B.
\end{cases}
\]
By~\eqref{eqBound3}, we see that $\by^B_{\ba} \ge (3-\gamma^{-1})\delta$ for all $\ba \in B$.
Thus, $\by^B \in \Z^n_{\ge 0}$.
By Lemma~\ref{lemLPBasisConstant}, $\by^B$ is optimal for $\LP(\bb - A \by^B)$.
The vector $\by^* - \by $ is also feasible for $\LP(\bb - A \by^B)$, so $\bc^\intercal (\by^* - \by) \le \bc^\intercal \by^B$.
The inequality $\by^B + \bu \ge \mathbf{0}^n$ holds because
\[
{\by}^{B}_{\ba} + \bu_{\ba} \ge (1-\gamma^{-1})\delta \ge 0 \qquad \forall ~ \ba \in B.
\]
This implies that $\by^B+\bu$ is feasible for $\LP(\bb - A \bw)$.
By Lemma~\ref{lemLPBasisConstant}, $\bw^B$ is optimal for $\LP(\bb - A \bw)$.
Therefore, $\bc^\intercal (\by^B + \bu) \le \bc^\intercal \bw^B$.
This proves~\eqref{eqFINALFINALEq1}.

It remains to prove~\eqref{eqFINALFINALEq}.
As $\by^B + \bu \ge \mathbf{0}^n$ and $\bw \ge \mathbf{0}^n$, it follows that
\[
({\by}^B+{\by})+ ({\bw} - {\by} + \bu) = ({\by}^{B}+\bu) + {\bw}
\]
is also nonnegative and feasible for $\IP(\bb)$. 
Note that 
\[
\bc^\intercal(\by^B+\by) \le \bc^\intercal \by^* = \bc^\intercal \by + \bc^\intercal (\by^*-\by)\le \bc^\intercal(\by^B+\by).
\]
Thus, $\by^B+\by$ is an optimal solution to $\IP(\bb)$ and $\bc^\intercal ({\bw} - {\by} + \bu) \le 0$. 
Because $\bz \in N(B)$, there exists ${\bb}^{\bz} \in D(B)$ and ${\bz}^B \in \Z^n_{\ge 0}$ such that $\{\ba \in A : {\bz}^B_{\ba} > 0\} \subseteq B$ and ${\bz}^B + \bz$ is optimal for $\IP({\bb}^{\bz})$. 
By~\eqref{eqBound3}, $\bz^B_{\ba} \ge (3-\gamma^{-1})\delta$ for all $\ba \in B$.
Hence, $\bz^B - \bu \ge \mathbf{0}^n$.
Recall that $\by \ge \mathbf{0}^n$, $\bz-\bw\ge \mathbf{0}^n$, and $A (\bw - \by  +\bu) = \mathbf{0}^m$ by definition.
Thus, 
\[
(\bz^B + {\bz}) - ({\bw} - {\by} + \bu)  = ({\bz}^{B} - \bu) +  (\bz - \bw) + {\by} 
\]
is feasible for $\IP({\bb}^{\bz}).$
This implies that $\bc^\intercal  ({\bw} - {\by} + \bu) \ge 0$.
\endproof
The final lemma in this section shows that certain vectors in $N(B)$ satisfy additional properties that we will use to prove Theorem~\ref{thmSuppProb}.
We notify the reader that the proof of Lemma~\ref{lemNonbasic3} is similar to the proof of Lemma~\ref{lemNonbasic2} although the main assumptions are different.

 \begin{lemma}\label{lemNonbasic3}
Let $B \subseteq A$ be an optimal $\LP$ basis matrix and $\bb \in D(B)$. 
Assume that $\bz$ minimizes $|\supp(\bz)|$ over all $\bz \in N(B)$ such that $A \bz \equiv_\Gamma \bb$.
If $\bw$ and $\by $ are distinct vectors satisfying $\bw_{\ba}, \by_{\ba} \in \{0, \bz_{\ba}\}$ for each $\ba \in A$, then $A \bw \not\equiv_\Gamma A \by$.
 \end{lemma}

\proof
Assume to the contrary that there exist distinct vectors $\bw$ and $\by$ such that $A \bw \equiv_\Gamma A \by$ and $\bw_{\ba}, \by_{\ba} \in \{0, \bz_{\ba}\}$ for each $\ba \in A$.
We may assume that $\supp({\bw}) \cap \supp({\by}) = \emptyset$ by subtracting the vector of overlapping support.
We assume without loss of generality that $\bw \neq \mathbf{0}^n$.
Note that $\bz - \bw + \by \in \Z^n_{\ge 0}$, $A(\bz - \bw + \by) \equiv_\Gamma \bb$, and $\supp(\bz - \bw + \by)$ is a strict subset of $\supp(\bz)$.
We cannot apply Lemma~\ref{lemNonbasic2} to conclude $\bz - \bw + \by \in N(B)$, which would contradict that $\bz$ had minimal support, because $\bz - \bw + \by \not \le \bz$.
Instead, we show that there exists a vector $\bv \in N(B)$ satisfying $\bv \le \bz - \bw + \by$ and $A \bv \equiv_{\Gamma} \bb$; this will yield the same contradiction.

Let $\bv \in \Z^n$ minimize $\|\bv\|_1$ over the integral vectors such that $\mathbf{0}^n \le \bv \le \bz - \bw + \by$ and $A \bv \equiv_{\Gamma} \bb $.
Condition \emph{(iii)} in~\eqref{eqDSet} is satisfied by $\bv$; otherwise, $\|\bv\|_1$ would not be minimized. 
To show that Conditions \emph{(i)} and \emph{(ii)} in~\eqref{eqDSet} hold, we define a suitable vector $\bv^B$. 
Define $\bv^B \in \R^n$  to be
\begin{equation}\label{eqFinalCrazyEq}
\bv^{B}_{\ba} := \begin{cases} 
 \left[B^{-1}( \bb - A \bv )\right]_{\ba} & \text{if}~ \ba \in B\\[.1 cm]
0 &\text{if}~ \ba \in A\setminus B.
\end{cases}
\end{equation}
By~\eqref{eqBound3} in Lemma~\ref{lemNonbasic1}, we have $\bv^B \in \Z^n_{\ge 0}$.
Also, $\{\ba \in A : \bv^B_{\ba} >0\} \subseteq B$ by construction. 
Hence, Condition \emph{(i)} in~\eqref{eqDSet} holds.

It is left to show Condition~\emph{(ii)} in~\eqref{eqDSet} holds, i.e., that $\bv^B+\bv$ is an optimal solution to $\IP(\bb)$.
By using the definition of $\bv^B$, it follows that $\bv^B+ \bv$ is feasible for $\IP(\bb)$.
It remains to show that $\bv^B+\bv$ is optimal.
Lemma~\ref{lemNonbasic2} applied to $\bz$ and $\bb$ implies that there exists a vector $\bz^B \in \Z^n_{\ge 0}$ such that $\{\ba \in A : \bz^B_{\ba} > 0\} \subseteq B$ and $\bz^B+\bz$ is optimal for $\IP(\bb)$. 
Because $A \bw  \equiv_{\Gamma} A \by$, there exists $\bu \in \Z^n$ such that
\[
\{\ba \in A : \bu_{\ba} \neq 0\} \subseteq B \quad \text{and} \quad A (\bw - \by + \bu) = \mathbf{0}^m.
\]
The argument used to prove~\eqref{eqFINALFINALEq} in the proof Lemma~\ref{lemNonbasic2} can be repeated to conclude
\(
\bc^\intercal(\bw-\by+\bu)  = 0 
\).
Hence, 
\[
\bc^\intercal (\bz^B+\bz) =  \bc^\intercal (\bz^B+\bz) - \bc^\intercal(\bw-\by+\bu) = \bc^\intercal \bv + \bc^\intercal [({\bz}^{B} - \bu) +  (\bz - \bw + {\by}) - \bv].
\]
If we can prove that
\begin{equation}\label{eqFinalSillBound}
\bc^\intercal [({\bz}^{B} - \bu) +  (\bz - \bw + {\by}) - \bv] \le  \bc^\intercal \bv^B ,
\end{equation}
then we will complete the proof that $\bv+\bv^B$ is optimal because
\[
\bc^\intercal (\bz^B+\bz)  \le \bc^\intercal (\bv^B+\bv).
\]
By~\eqref{eqBound3}, $\bz^B_{\ba} \ge (3-\gamma^{-1})\delta$ for each $\ba \in B$. 
Using the facts that $\bw$ and $\by$ have disjoint supports and that $\bw_{\ba}, \by_{\ba} \in \{0, \bz_{\ba}\}$ for each $\ba \in A$, we have $\|\bw - \by\|_1 \le \|\bz\|_1$.
Thus,
\[
\|\bu\|_{\infty} = \| B^{-1}A(\bw - \by)\|_{\infty} \le \|B^{-1} A\|_{\infty} \cdot \|\bw - \by\|_1 
\le \|B^{-1} A\|_{\infty} \cdot \|\bz\|_1  \le \gamma^{-1} \cdot \delta
\]
and $\bz^B - \bu \ge \mathbf{0}^n$. 
Moreover, $({\bz}^{B} - \bu) +  (\bz - \bw + {\by}) - \bv \ge \mathbf{0}^n $ because $\mathbf{0}^n \le \bv \le \bz - \bw + {\by}$. 
Finally, $({\bz}^{B} - \bu) +  (\bz - \bw + {\by}) - \bv $ and $\bv^B$ are both feasible for $\LP(A \bv^B)$ with $\bv^B$ being optimal by Lemma~\ref{lemLPBasisConstant}. 
This proves~\eqref{eqFinalSillBound}.
\endproof

\section{Results about $\sigma$ and $\pi$}\label{secSupp}
%

Our remaining goal is to complete the proofs of Theorem~\ref{thmSuppProb} and Theorem~\ref{thmMainProx}.
We proceed as follows in both proofs.
Define $\Lambda := A \cdot \Z^m$. 
Let $A^1, \dotsc, A^s \subseteq A$ be the optimal $\LP$ basis matrices.
By Lemma~\ref{lemLPBasisConstant}, these matrices form a simplicial covering of $\cone(A)$.
As in~\eqref{eqGammaLattice},~\eqref{eqBSet}, and~\eqref{eqDSet}, define
\[
\Gamma^i := A^i \cdot \Z^m, ~~D^i := D(A^i), ~~\text{and}~~ N^i := N(A^i) \qquad \forall~ i \in \{1, \dotsc, s\}.
\]
In view of \eqref{eqBSet}, we define the vectors $\bd^i :=A^i(3 \delta\cdot \mathbf{1}^m)$ for all $i\in \{1,\dotsc,s\}$.

\proof[Proof of Theorem~\ref{thmSuppProb}]
In accordance with equation \eqref{thmMain:condition} from Theorem~\ref{thmMain}, we define the set
\[
X^i:=\left\{\bg \in \Lambda/\Gamma^i : 
\max\left\{\sigma(\bb) :
\begin{array}{l}\bb \equiv_{\Gamma^i} \bg,\\[.05 cm]
 \bb\in \cone(A^i) + \bd^i
 \end{array}
\right\} \le m+k \right\}
\]
and show that
\begin{equation}\label{eqSupportInduction1}
| X^i | \ge  \min\left\{\gamma^{-1} \cdot |\det(A^i)|, 2^k \right\} \qquad \forall ~i \in \{1, \dotsc, s\}.
\end{equation}
Theorem~\ref{thmSuppProb} then follows from Theorem~\ref{thmMain} with $\alpha = m+k$ and $\beta_i \ge \min\{\gamma^{-1} \cdot |\det(A^i)|,$ $ 2^k \}$.

Fix $i \in \{1, \dotsc, s\}$.
We complete the proof of~\eqref{eqSupportInduction1} in two cases.

\smallskip

\noindent \textbf{Case 1.}
Assume that $\Lambda / \Gamma^i = X^i$. 
By~\eqref{eqNormalizedGCD}, we have
\[
|X^i|  = |\Lambda / \Gamma^i| = \gamma^{-1} \cdot |\det(A^i)|.
\]
This proves~\eqref{eqSupportInduction1}.

\smallskip

\noindent \textbf{Case 2.}
Assume that $\Lambda / \Gamma^i \supsetneq X^i$. 
By the definition of $X^i$, there exists $\bg \in \Lambda / \Gamma^i$ such that
\[
\max\left\{\sigma(\bb) :
\hspace{-.15 cm}
\begin{array}{l}\bb \equiv_{\Gamma^i} \bg,\\[.05 cm]
 \bb\in \cone(A^i) + \bd^i
 \end{array} \right\}> m+k.
\]
Lemma~\ref{lemNonbasic2} implies that for any $\bb \in \cone(A^i) + \bd^i$ and any $\bz^{\bg} \in N^i$ with $A \bz^{\bg} \equiv_{\Gamma^i} \bg$, there exists an optimal solution to $\IP(\bb)$ whose support is bounded by $m + |\supp( \bz^{\bg} )|$. 
Hence, 
\begin{equation}\label{eqTooBigSupp}
\min\left\{|\supp(\bz^{\bg})| : \bz^{\bg} \in N^i ~\text{and}~ A\bz^{\bg} \equiv_{\Gamma^i} \bg\right\} \ge k+1.
\end{equation}
Choose $\bg$ and $\bz^{\bg} \in N^i$ as argument maximizers and minimizers, respectively, of the problem
\[
\max_{\bg \in \Lambda / \Gamma^i} \min\left\{|\supp(\bz^{\bg})| : \bz^{\bg} \in N^i ~\text{and}~ A\bz^{\bg} \equiv_{\Gamma^i} \bg\right\}.
\]
Inequality~\eqref{eqTooBigSupp} implies that $|\supp(\bz^{\bg})| \ge k+1$.

Define the sets
\[
Z^i:=\{\bz\in\Z^n : \bz_{\ba} \in \{0, \bz^{\bg}_{\ba}\}\text{ for each }\ba \in A  \text{ and } |\supp({\bz})| \le k\}
\]
and
\[
H^i := \{ \bh \in \Lambda / \Gamma^i : \bh \equiv_{\Gamma^i} A \bz \text{ for some  } {\bz} \in Z^i\}.
\]
We show that $H^i \subseteq X^i$.
Let $\bh \in H^i$ and take $\bb \in \cone(A^i) + \bd^i$ such that $\bb \equiv_{\Gamma^i} \bh$.
There exists $ {\bz} \in Z^i$ such that $A {\bz} \equiv_{\Gamma^i} \bb$. 
The definition of $N^i$ and Lemma~\ref{lemNonbasic2} imply that there exists an optimal solution to $\IP(\bb)$ of the form $\bz + \bz^{i}$, where
\(
\{\ba \in A: \bz^{ i}_{\ba} > 0 \} \subseteq A^i.
\)
Hence, 
\[
\sigma(\bb) \le |\supp(\bz+\bz^{i} )| \le |\supp(\bz^{i})| +|\supp(\bz)| \le m+k.
\]
This implies that $H^i \subseteq X^i$.
As $\bz^{\bg}$ was chosen to have minimal support, it follows from Lemma~\ref{lemNonbasic3} that $Z^i$ and $H^i$ have the same cardinality.
Thus,
\begin{equation}\label{eq:for:example}
|X^i| \ge |H^i| = |Z^i| =   \sum_{j=0}^k {|\supp(\bz^{\bg})| \choose j} \ge \sum_{j=0}^k {k+1 \choose j} \ge \sum_{j=0}^k {k \choose j}= 2^k.
\end{equation}
\endproof

\begin{remark}\label{remSparsFeas}
If $\bc = \mathbf{0}^n$, then $\sigma(\bb)$ denotes the sparsest feasible solution to $\IP(\bb)$.
Under this assumption, every invertible matrix $B \subseteq A$ is an optimal $\LP$ basis matrix, and we can upper bound asymptotic densities of $\sigma$ in terms of the \emph{smallest positive determinant} of all the submatrices of $A$.
Define
$$
\eta :=  \min ~  \{|\det(B)|: B \subseteq A \text{ is invertible}\},
$$
and let $B \subseteq A$ be a matrix that attains this minimum.
Suppose $A^1, \dotsc, A^s \subseteq A$ form a simplicial covering of $\cone(A)$.
Provided $\bb$ is deep in $\cone(A^i)$, one can express $\bb$ as $\bb=A^i\bz+B\by$, where $\bz\in\Z^m_{\ge0}$ and $\by\in\R^m_{\ge0}$.
Following the proof of Theorem~\ref{thmSuppProb}, for every fixed vector $\bz \in \Z^m$, it holds that
\[
\Pr\left(\{ \bb \in A^i \bz +\cone(B) : \sigma(\bb) \le 2m+k\} \right) \ge \frac{2^k}{\gamma^{-1} \cdot \eta}.
\]
The term $2m + k$ comes from two places: $m+k$ is from Theorem~\ref{thmSuppProb}, and the extra $m$ comes from $\bz \in \Z^m_{\ge0}$. 
Because this bound holds for every $\bz \in \Z^m_{\ge0}$ and the basis matrix $A^i$ was arbitrarily chosen, we can let $\bz$ vary to cover the deep regions corresponding to every basis matrix. 
Thus,
\[
\Pr\left( \sigma \le 2 m + k \right) 
\ge  \min\bigg\{1, ~\frac{2^k}{\gamma^{-1} \cdot \eta}\bigg\}.
\]
This is closely related to the results on the sparsity of systems of linear Diophantine equations in~\cite{AlAvDeOe19}.
\end{remark}

\medskip

\proof[Proof of Theorem~\ref{thmMainProx}]

We first prove Part~{\it(a)}.
In accordance with \eqref{thmMain:condition} from Theorem~\ref{thmMain}, we define the set
\begin{align*}
X^i&:=\left\{\bg \in \Lambda/\Gamma^i : 
\max\left\{\pi(\bb) :
\begin{array}{l}\bb \equiv_{\Gamma^i} \bg,\\[.05 cm]
 \bb\in \cone(A^i) + \bd^i
 \end{array}
\right\} \le m \gamma^{-1}\cdot \delta\cdot \frac{k}{k+1} + k\right\}
\end{align*}
and show that
\begin{equation}\label{eqProxInduction1}
| X^i | \ge  \min\left\{\gamma^{-1}\cdot |\det(A^i)|, ~k+1\right\} \qquad \forall ~i \in \{1, \dotsc, s\}.
\end{equation}
The result then follows from Theorem~\ref{thmMain}.
Fix $i \in \{1, \dotsc, s\}$.

\smallskip
\noindent\textbf{Case 1.} Assume that $ \Lambda / \Gamma^i = X^i$.
By~\eqref{eqNormalizedGCD}, we have
\[
|X^i| = \gamma^{-1}\cdot |\det(A^i)|.
\]
This shows~\eqref{eqProxInduction1}.

\smallskip
\noindent\textbf{Case 2.} 
Assume that $ \Lambda / \Gamma^i \supsetneq X^i$. 
Consider any $\bg \in \Lambda / \Gamma^i$, $\bz^{\bg} \in N^i$, and $\bb \in \cone(A^i) + \bd^i$ such that $\bg  \equiv_{\Gamma^i} A\bz^{\bg} \equiv_{\Gamma^i} \bb$.
Lemma~\ref{lemNonbasic2} implies that there exists an optimal solution to $\IP(\bb)$ of the form $\bz^{\bg}+\bz^{i}$, where
\(
\{\ba \in A : \bz^{i}_{\ba} > 0\} \subseteq A^i.
\)
Let $\bx^*$ be the optimal vertex solution to the linear program $\LP(\bb)$ with $\{ \ba \in A : \bx^*_{\ba} > 0\} \subseteq A^i$.
The supports of $\bx^*$ and $\bz^i$ are contained in $A^i$ while the support of $\bz^{\bg} \in N^i$ is disjoint from $A^i$ by Condition~\emph{(iii)} in~\eqref{eqDSet}. 
Hence, the supports of $\bx^*-\bz^i$ and $\bz^{\bg}$ are disjoint. 
From this and~\eqref{eqBound2}, we see that
\begin{equation}\label{eqProxBoundZg}
\begin{aligned}
\pi(\bb)  = \|\bx^* - (\bz^{i}+\bz^{\bg})\|_{1}  &=  \|\bx^* -\bz^{i}\|_{1}+\|\bz^{\bg}\|_{1} = \|(A^i)^{-1}A\bz^{\bg}\|_{1}+\|\bz^{\bg}\|_{1}\\[.15 cm]
& \le m\cdot \frac{\delta}{|\det(A^i)|}\cdot \|\bz^{\bg}\|_{1}+\|\bz^{\bg}\|_{1}.
\end{aligned}
\end{equation}
Because $\Lambda / \Gamma^i \supsetneq X^i$, there exists a particular $\bg \in \Lambda / \Gamma^i$ such that
\begin{align*}
\max\left\{\pi(\bb) :
\begin{array}{l}\bb \equiv_{\Gamma^i} \bg,\\[.05 cm]
 \bb\in \cone(A^i) + \bd^i
 \end{array}
\right\} &> m \gamma^{-1}\cdot \delta\cdot \frac{k}{k+1} + k.
%
\intertext{Let $\bz^{\bg} \in N^i$ satisfy $A\bz^{\bg} \equiv_{\Gamma^i} \bg$. 
By~\eqref{eqProxBoundZg}, we have}
%
\max\left\{\pi(\bb) :
\begin{array}{l}\bb \equiv_{\Gamma^i} \bg,\\[.05 cm]
 \bb\in \cone(A^i) + \bd^i
 \end{array}\right\} &\le m \cdot \frac{\delta}{|\det(A^i)|} \cdot \|\bz^{\bg}\|_1 + \|\bz^{\bg}\|_1.
\end{align*}
If $\|\bz^{\bg}\|_1 < k$, then the latter two inequalities imply that 
\[
m \cdot \frac{\delta}{|\det(A^i)|} \cdot \|\bz^{\bg}\|_1 + \|\bz^{\bg}\|_1
 > m \gamma^{-1}\cdot \delta\cdot \frac{k}{k+1} + k
> m \gamma^{-1}\cdot \delta\cdot \frac{\|\bz^{\bg}\|_1}{\|\bz^{\bg}\|_1+1} + \|\bz^{\bg}\|_1,
\]
or equivalently that $\|\bz^{\bg}\|_1 \ge \gamma^{-1} \cdot |\det(A^i)|$. 
However, this contradicts~\eqref{eqBound1}.
Hence, $\|\bz^{\bg}\|_1 \ge k$ and $\gamma^{-1} \cdot |\det(A^i)| \ge k+1$.

Let $\bz \in \Z^n$ satisfy $\mathbf{0}^n \le \bz \le \bz^{\bg}$ and $\|\bz\|_1 = k$. 
Consider the set
\[
H^i := \{ \bh \in \Lambda / \Gamma^i : \bh \equiv_{\Gamma^i} A \overline{\bz} \text{ for }\overline{\bz} \in \Z^n \text { with } \mathbf{0}^n \le \overline{\bz} \le\bz \}.
\]
We claim that $H^i \subseteq X^i$.
Take $\bh \in H^i$ and let $\overline{\bz} \in \Z^n$ satisfy $\mathbf{0}^n \le \overline{\bz} \le\bz$ and $\bh \equiv_{\Gamma^i} A\overline{\bz}$.
By Lemma~\ref{lemNonbasic2}, both $\bz$ and $\overline{\bz}$ are in $N^i$.
Let $\bb \in \cone(A^i) + \bd^i$ be such that $\bb \equiv_{\Gamma^i} \bh$.
Applying~\eqref{eqProxBoundZg} to $\overline{\bz}$, it follows that 
\begin{align*}
\pi(\bb) 
= \|(A^i)^{-1}A\overline{\bz}\|_1 + \|\overline{\bz}\|_1 
&\le m \cdot \frac{ \delta}{|\det(A^i)|} \cdot \|\overline{\bz}\|_1 + \|\overline{\bz}\|_1\\
 &\le  m \cdot \gamma^{-1} \cdot \delta \cdot \frac{ k}{k+1}+k.
\end{align*}
Hence, $\bh \in X^i$ and $H^i \subseteq X^i$.
Because $\bz \in N^i$, Condition \emph{(iii)} in~\eqref{eqDSet} implies that $A \bv \not\equiv_{\Gamma^i} A \bw$ for every $\bv,\bw \in \Z^n$ satisfying $\mathbf{0}^n\le \bv\lneq\bw\le\bz$.
Therefore,
\[
|X^i| \ge 
 |H^i| \ge \|\bz\|_1 + 1 = k+1 \ge  \min\{\gamma^{-1}\cdot |\det(A^i)|,~ k+1\},
\]
which completes the proof of~\eqref{eqProxInduction1} and proves Part \emph{(a)} of the theorem.

The proof of Part~{\it(b)} is almost identical to the proof of Part~{\it(a)}.
One defines
\begin{align*}
X^i_{\infty}&:=\left\{\bg \in \Lambda/\Gamma^i : 
\max\left\{\pi^{\infty}(\bb) :
\begin{array}{l}\bb \equiv_{\Gamma^i} \bg,\\[.05 cm]
 \bb\in \cone(A^i) + \bd^i
 \end{array}
\right\} \le \gamma^{-1}\cdot \delta\cdot \frac{k}{k+1} \right\}
\end{align*}
and shows that
\[
| X^i_{\infty} | \ge  \min\left\{\gamma^{-1}\cdot |\det(A^i)|, ~k+1\right\} \qquad \forall ~i \in \{1, \dotsc, s\}.
\]
The key difference is that we replace \eqref{eqProxBoundZg} with
\begin{align*}
\pi^{\infty}(\bb)  =  \max\left\{ \|\bx^* - \bz^i\|_{\infty}, \|\bz^{\bg}\|_{\infty}\right\}
&= \max\left\{ \|(A^i)^{-1}A\bz^{\bg}\|_{\infty}, \|\bz^{\bg}\|_{\infty}\right\}\\[.1 cm]
& \le \max\left\{ \|(A^i)^{-1}A\|_{\infty}\|\bz^{\bg}\|_{1}, \|\bz^{\bg}\|_{\infty}\right\}\\[.1 cm]
& = \|(A^i)^{-1}A\|_{\infty}\|\bz^{\bg}\|_{1}\\[.1 cm]
& \le  \frac{\delta}{|\det(A^i)|} \cdot \|\bz^{\bg}\|_{1} .
\end{align*}
\endproof

\smallskip
\begin{remark}\label{remark:uniqueOptima}
In Section~\ref{subsecProx}, we made the assumption that the optimal solution to $\LP(\bb)$ is unique for all feasible $\bb$. 
If this assumption is dropped, then the definition of distance should be adapted as follows.
Define the \emph{minimum distance between an optimal $\LP$ vertex solution and an optimal $\IP$ solution} to be
\[
 \pi^{\min}(\bb)
 :=  \min_{\bx^*} ~ \min_{\bz^*}\bigg\{\|\bx^* - \bz^*\|_1 : 
\begin{array}{l}
\bx^* \text{ is an optimal vertex solution to } \LP(\bb)\\
\bz^* \text{ is an optimal solution to } \IP(\bb)
\end{array}
\bigg\},
\]
and the \emph{maximum of the minimum distance between $\LP$ optimal vertices and $\IP$ optimal solutions} to be
\[
 \pi^{\max}(\bb)
 :=  \max_{\bx^*} ~ \min_{\bz^*}\bigg\{\|\bx^* - \bz^*\|_1 : 
\begin{array}{l}
\bx^* \text{ is an optimal vertex solution to } \LP(\bb)\\
\bz^* \text{ is an optimal solution to } \IP(\bb)
\end{array}
\bigg\}.
\]
If $\IP(\bb)$ is infeasible, then $\pi^{\min}(\bb)  = \pi^{\max}(\bb):= \infty$.
The value $\pi^{\min}(\bb)$ can be bounded by considering only one solution to $\LP(\bb)$ while $\pi^{\max}(\bb)$ needs to consider every optimal vertex of $\LP(\bb)$.
It follows immediately from Theorem~\ref{thmMainProx} that
\[
\Pr\left(\pi^{\min} \le m \gamma^{-1} \cdot \delta\cdot \frac{k}{k+1} + k \right)  \ge  \frac{k+1}{\gamma^{-1} \cdot \delta} \qquad \forall ~ k \in \{0, \dotsc, \gamma^{-1} \cdot \delta - 1\}.
\]
It is not clear if $\pi^{\max}(\bb)$ can be bounded in the same way.
However, for the extreme case $k = \gamma^{-1} \cdot \delta - 1$ it can be shown that
\[
\Pr\left(\pi^{\max} \le  (m+1) (\gamma^{-1} \cdot \delta-1)\right) =1.
\]
The proof of this equation is similar to the proof of Theorem~\ref{thmMainProx}, and it is omitted here.
\end{remark}

\begin{remark}
As a final remark, we want to point out that our proofs provide a method for computing exact densities.
Let us illustrate this by considering again the sparsity function $\sigma$.
Set $\bc=(\mathbf{1}^m,\mathbf{0}^m)$ and $A=[2I,I]$, where $I$ denotes $m\times m $ identity matrix. 
There is a unique optimal $\LP$ basis matrix, which is defined by the first $m$ columns.
The asymptotic densities for $k=0,1\ldots,m$ are
\[
\Pr(\sigma\le m+k)=\frac{1}{2^m}\sum_{i=0}^k {m \choose i},
\]
which can be inferred from \eqref{eq:for:example}.
Note that this coincides with Theorem~\ref{thmSuppProb} for $k=0$.
\end{remark}

\section*{Acknowledgments}
The authors would like to thank Laurence Wolsey, Luze Xu, and the anonymous referees for helping us greatly improve the presentation of the material.
The third author was supported by the Einstein Foundation Berlin.

\bibliographystyle{siamplain}
\bibliography{references}
\end{document}